\documentclass[reqno,a4paper,11pt]{article}

\usepackage{epic,eepic,epsfig,amssymb,amsmath,amsfonts,amsthm,mathrsfs}

\title{Mass transportation with LQ cost functions}

\author{A.~Hindawi\footnote{Universit\'e de Nice-Sophia
  Antipolis, Labo.\ J.-A.\ Dieudonn\'e, UMR 6621, Parc
  Valrose, 06108 Nice Cedex 02, France \textit{and} INRIA Sophia Antipolis. ({\tt
    hindawi@unice.fr})}
\and
L.~Rifford\footnote{Universit\'e de Nice-Sophia
  Antipolis, Labo.\ J.-A.\ Dieudonn\'e, UMR 6621, Parc
  Valrose, 06108 Nice Cedex 02, France. On leave to INRIA Sophia Antipolis for one year. ({\tt
    rifford@unice.fr})}
\and
J.-B.~Pomet\footnote{INRIA, B.P. 93, 06902 Sophia
  Antipolis cedex, France. ({\tt Jean-Baptiste.Pomet@sophia.inria.fr
   })}}

\date{}
\numberwithin{equation}{section}

\newtheorem{theorem}{Theorem}[section]

\newtheorem{proposition}[theorem]{Proposition}

\theoremstyle{definition}
\theoremstyle{definition}\newtheorem{definition}[theorem]{Definition}

\theoremstyle{definition}

\def\R{\textrm{I\kern-0.21emR}}
\def\N{\textrm{I\kern-0.21emN}}

\def\p{\partial}
\def\g{\gamma}

\newcommand{\supp}{{\rm \operatorname{supp}}}

\renewcommand{\exp}{{\rm exp}}

\begin{document}

\maketitle

\begin{abstract}
We study the optimal transport problem in the Euclidean space where the cost function is given by the value function associated with a Linear Quadratic minimization problem. Under appropriate assumptions, we generalize Brenier's Theorem proving existence and uniqueness of an optimal transport map. In the controllable case, we show that the optimal transport map has to be the gradient of a convex function up to a linear change of coordinates. We give regularity results and also investigate the non-controllable case.
\end{abstract}

\section{Introduction}

The optimal transport problem can be stated as follows: given two probability measures $\mu_0$ and $\mu_1$, defined on measurable
spaces $X$ and $Y$ respectively and a cost function 
$$
c \, : \, X \times Y \, \longrightarrow \, \R \cup \{+ \infty\}
$$
find a measurable map 
$$
T \, : \, X \, \longrightarrow \, Y
$$
which pushes forward $\mu_0$ to $\mu_1$, that is 
$$
T_\sharp \mu_0=\mu_1 \qquad \text{(i.e. $\mu_1(B)=\mu_0\left(T^{-1}(B)\right)$ for all $B \subset Y$ measurable)},
$$
and which minimizes the transportation cost
$$
\mbox{cost}_c (T) := \int_X c(x,T(x))\,d\mu_0(x).
$$
When the transport condition $T_\sharp \mu_0=\mu_1$ is satisfied, we say that $T$ is a \textit{transport map}, and if $T$ minimizes also the cost, that is if
$$
\mbox{cost}_c (T)  = \min_{S_\sharp \mu_0=\mu_1} \Bigl\{ \mbox{cost}_c (S) \Bigr\}
$$
then we call it an \textit{optimal transport map}. Since the seminal famous paper by Gaspard Monge in 1781 \cite{monge}, there was a revival of interest in mass transportation in the nineties.  In 1987 \cite{brenier87,brenier91}, Brenier proved an existence and uniqueness result of optimal transport maps for the cost $c(x,y)= |x-y|^2$ in $\R^n$ and showed that any such optimal transport map is indeed the gradient of a convex function. Since then, people extended the theory to other costs functions in $\R^n$ or to other types of spaces (see \cite{villaniSF}). 

The aim of this paper is to study existence, uniqueness, and regularity of optimal transport maps for costs functions coming for LQ minimization problems in $\R^n$. 
Let us consider a linear control system of the form
\begin{eqnarray}\label{system}
\dot{x}= A x + B u
\end{eqnarray}
where the state $x$ belongs to $\R^n$, the control $u$ belongs to $\R^m$ and  $A,B$ are $n \times n $ and  $n \times m$ matrices respectively. For every initial state $x\in \R^n$ and every control $u \in L^2\bigl( [0,1];\R^m\bigr)$, we denote by $x(\cdot;x,u) : [0,1] \rightarrow \R^n$ the unique solution to the Cauchy problem
\begin{eqnarray}\label{cauchy}
\left\{ 
\begin{array}{l}
\dot{x}(t) = A x(t) + B u(t) \quad \mbox{for a.e. } t \in [0,1],\\
x(0) =x.
\end{array}
\right.
\end{eqnarray}
Let us in addition consider a quadratic Lagrangian of the form 
\begin{eqnarray}\label{lagrangian}
L(x,u)= \frac{1}{2} \langle x, W x\rangle  + \frac{1}{2} \langle u, U u \rangle,  
\end{eqnarray} 
where $W$ is a $n \times n$ symmetric non-negative matrix and $U$ is a $m \times m$ symmetric positive definite matrix. The cost $c : \R^n \times \R^n \rightarrow [0,+\infty]$ associated with (\ref{system}) and (\ref{lagrangian}) is given by 
\begin{eqnarray}\label{costLQ}
c(x,y) := \inf \left\{  \int_{0}^{1} L( x(t;x,u),u(t)) dt  \, \vert \, u\in L^2\bigl( [0,1];\R^m\bigr) \, \mbox{ s.t. } x(1;x,u) = y \right\},
\end{eqnarray}
where we set $c(x,y)=+\infty$ if there is no $u\in L^2\bigl( [0,1];\R^m\bigr)$ such that $x(1;x,u)=y$.

We notice that LQ costs as above include as a particular case the Euclidean cost $1/2 |x-y|^2$ by taking $B=I_n$, $U=I_n$, $A=0$, $W=0$ . 
Costs coming from optimal control have already been studied in \cite{al07} (see also \cite{fr10}).  However, this reference does not give regularity results or properties relating the transport map with the gradient of a convex function. Furthermore, the study of a cost that is finite only for points lying in the same leaf of a foliation, like in the non-controllable case is also original. 

The structure of the paper is the following. Our results are stated in Section \ref{RESULTS}. In Section \ref{SPRE}, we introduce some preliminaries in optimal transport theory concerned with Kantorovitch duality. In Sections \ref{SPROOFTHM1} and \ref{SPROOFTHM2}, we provide the proofs of our results. Then, we present  some examples in Section \ref{SEXAMPLES}. Finally in Section \ref{Scomments}, we conclude with several remarks about our results. 

\section{Main results}
\label{RESULTS}

\subsection{Preliminaries on linear systems}

Consider system (\ref{system}) and let $V$ be the smallest linear subspace of $\R^n$ that contains the image
of the operator $B:\R^m\to\R^n$ and is invariant by $A$, or, more explicitly,
\begin{equation}
\label{eq:kalman}
V=\mathrm{Span}_\R \Bigl\{B,AB,A^2B,\ldots,A^{n-1}B \Bigr\}   \subset\R^n\,,\ \ d=\dim V.
\end{equation}
The well-known \emph{Kalman criterion} states that system (\ref{system}) is
controllable if and only if $d$; this is stated below with a description of the situation where $d<n$.
This decomposition can be found for instance in \cite[Lemma 3.3.3]{sontag90}.
\begin{proposition}
 \label{prop-cont}
 Let $x,y$ be in $\R^n$. 
 There exists a control $u \in L^2\bigl([0,1];\R^m\bigr)$ such that $x(t;x,u)=y$ if and only if $e^{A}x-y\in
 V$ (i.e. $y$ lies in the affine subspace $e^{A}x+V$). In particular it exists for all $x,y$ if $d$; the system --or the pair $(A,B)$-- is then called \emph{controllable}.


 If $d<n$, after a linear change of coordinates in $\R^n$, the control system (\ref{system}) has the
 following form:
\begin{eqnarray}\label{system2}
\dot{x} = \left( \begin{matrix}
\dot{x}_1\\
\dot{x}_2 
\end{matrix}
\right) = \left( \begin{matrix}
A_1 & A_3 \\ 0 & A_2 
\end{matrix}
\right)  \left( \begin{matrix}
x_1 \\ x_2
\end{matrix}
\right) + \left( \begin{matrix}
B_1 u\\
0
\end{matrix}
\right),
\end{eqnarray}
where the state $x$ is partitioned into two blocks $x_1$ and $x_2$ of dimension $d$ and $n-d$ respectively
and $A_1$, $A_2$, $A_3$, $B_1$ are $d\times d$, $d\times n-d$, $n-d\times n-d$ and $d\times m$ matrices
respectively, such that the pair $(A_1,B_1)$ is controllable.
\end{proposition}

\subsection{The controllable case}

The following result of existence, uniqueness and regularity follows easily from the classical theory of optimal mass transportation. Throughout the paper, $M^*$ denotes the transpose of the matrix $M$.

\begin{theorem}\label{THM1} 
Assume that the linear control system (\ref{system}) is controllable. Then there are symmetric positive
definite $n \times n$ matrices $D$ and $F$, and an invertible $n\times n$ matrix $E$ such that
\begin{eqnarray}\label{costTHM1}
c(x,y) =  \frac{1}{2} \left\langle x, D \, x \right\rangle -  \left\langle x, E \, y \right\rangle +  \frac{1}{2} \left\langle y, F \, y \right\rangle \qquad \forall x,y \in \R^n.
\end{eqnarray}
Let $\mu_0, \mu_1$ be two compactly supported probability measures on $\R^n$. Assume that $\mu_0$ is absolutely continuous with respect to the Lebesgue measure $\mathcal{L}^n$. Then, there is existence and uniqueness of an optimal transport map $T:\R^n \rightarrow \R^n$. Such a map is characterized by the existence of a convex function $\varphi: \R^n \rightarrow \R$ such that 
\begin{eqnarray}\label{TransportTHM1}
T(x) = E^{-1} \, \nabla \varphi (x) \qquad \mbox{for a.e. } x \in \R^n.
\end{eqnarray}
If in addition, $\mu_0, \mu_1$ are associated with probability densities $f_0, f_1$ on $\supp(\mu_0)$ and $\supp(\mu_1)$ respectively, with $f_0$ and $f_1$ bounded from below and above, and if  $\supp(\mu_0)$ is connected and $\supp(\mu_1)$ convex, then $T$ is continuous.
\end{theorem}

We postpone several remarks concerning Theorem \ref{THM1} to Section \ref{Scomments}.

\subsection{The non-controllable case}

Here $d<n$ in (\ref{eq:kalman}), hence the evolution of the $(n-d)$-dimensional $x_2$ is totally fixed by (\ref{system2}) and $c(x,y)$ will obviously be infinite for almost all pairs $x,y$: these such that $y_2\neq e^{A_2}x_2$.
A result similar to the above theorem requires first there exists at least one transport map with finite cost, i.e.
\begin{eqnarray}\label{HYPTHM2}
\inf_{S_\sharp \mu_0=\mu_1} \Bigl\{ \mbox{cost}_c (S) := \int_{\R^n} c(x,S(x))\,d\mu_0(x),    \Bigr\} < + \infty. 
\end{eqnarray}
Let $\pi_2:\R^n\to\R^{n-d}$ be the projection on the second block in the coordinates of (\ref{system2}): $\pi_2(x)=x_2$ (it is also the projection on the quotient $\R^n/V$).

\begin{theorem}\label{THM2} 
Let $\mu_0$ and $\mu_1$ be compactly supported probability measures with \emph{continuous} densities $f_0,
f_1$ ($\mu_0= f_0 \mathcal{L}^n, \mu_1= f_1 \mathcal{L}^n$).
There exists a transport map with finite cost (\ref{HYPTHM2}) if and only if
\begin{eqnarray}\label{HYPTHM2bis}
(\pi_2)_\sharp \mu_1=\left(e^{A_2}\circ\pi_2\right)_\sharp \mu_0\,.
\end{eqnarray}
If this is satisfied, then there exists a unique optimal transport map $T : \R^n \rightarrow \R^n$. Moreover, if both densities $f_0, f_1$ are bounded from below and above on $\supp(\mu_0)$ and $\supp(\mu_1)$ respectively, and if $\supp(\mu_0)$ and $\supp(\mu_1)$ are both convex, then $T$ is continuous.
\end{theorem}

Again, remarks concerning Theorem \ref{THM2} are postponed to Section \ref{Scomments}.\\

\section{Preliminaries in optimal transport theory}\label{SPRE}

Given two probability measures $\mu_0, \mu_1$ on $\R^n$ and a cost function $c: \R^n\times \R^n \to [0,+\infty]$, we are looking for a transport map
$T:\R^n\to \R^n$ which minimizes the transportation cost $\int_{\R^n} c(x,T(x))\,d\mu_0$. The constraint $T_\#\mu_0=\mu_1$  being highly non-linear, the optimal transport problem is quite difficult from the viewpoint of calculus
of variation. The major advance on this problem was due to
Kantorovich, who proposed in \cite{kant1, kant2} a notion of weak
solution of the optimal transport problem. He suggested to look
for \textit{plans} instead of transport maps, that is probability
measures $\gamma$ in $\R^n \times \R^n$ whose marginals are $\mu_0$ and
$\mu_1$, that is
$$
(\pi_1)_\sharp \gamma=\mu_0 \qquad \text{and} \qquad (\pi_2)_\sharp
\gamma=\mu_1,
$$
where $\pi_1: \R^n \times \R^n  \rightarrow \R^n$ are the canonical projections on the first and second variable respectively. Denoting by
$\Pi(\mu_0,\mu_1)$ the set of plans, the new minimization problem 
becomes the following:
\begin{equation}
\label{kantprob0} C \bigl(\mu_0,\mu_1\bigr) =\min_{\gamma \in \Pi(\mu_0,\mu_1)}
\left\{ \int_{\R^n \times \R^n} c(x,y)\,d\gamma(x,y) \right\}.
\end{equation}
If $\gamma$ is a minimizer for the Kantorovich formulation, we say
that it is an \textit{optimal plan}. Due to the linearity of the
constraint $\gamma \in \Pi(\mu_0,\mu_1)$, it is simple, using weak
topologies, to prove existence of solutions to (\ref{kantprob0}) 	as soon as $c$ is lower semi-continuous (see for instance \cite{villaniSF}). The connection between the formulation of Kantorovich
and that of Monge can be seen by noticing that any transport map
$T$ induces the plan defined by $(Id \times T)_\sharp \mu_0$ which
is concentrated on the graph of $T$. Thus, the problem of showing
existence of optimal transport maps can be reduced to prove that an
optimal transport plan is concentrated on a graph. Moreover, if one can show that \textit{any}
optimal plan is concentrated on a graph, since $\frac{\g_1+\g_2}{2}$ is optimal if so are $\g_1$ and $\g_2$,
uniqueness of the transport map easily follows. The following definition is exactly \cite[Definition 5.2]{villaniSF}:

\begin{definition}\label{c-convex}
{\rm
Let $c: \R^n\times \R^n \rightarrow [0,\infty]$. A function $\psi: \R^n \to \R \cup \{+\infty\}$ is said to be \textit{$c$-convex} if it is not identically $+\infty$ and there exists a function $\zeta :\R^n \rightarrow \R \cup \{-\infty,+\infty\}$ such that
$$
\psi(x) = \sup_{y\in \R^n} \Bigl( \zeta(y) - c(x,y) \Bigr) \qquad \forall x \in \R^n.
$$
Then its \textit{$c$-transform} is the function $\psi^c$ defined by
$$
\psi^c(y)=\inf_{x \in \R^n} \Bigl( \psi(x) + c(x,y) \Bigr) \qquad \forall y\in \R^n,
$$
and its \textit{$c$-subdifferential} is the set defined by
$$
\p_c\psi := \Bigl\{ (x,y) \in \R^n\times \R^n \, \vert \,  \psi^c(y)-\psi(x)=c(x,y) \Bigr\}.
$$
Moreover, the \textit{$c$-subdifferential of $\psi$} at $x\in \R^n$ is
$$
\p_c\psi (x) :=\Bigl\{ y\in  \R^n \, \vert \, (x,y) \in \p_c\psi \Bigr\},
$$
or equivalently 
$$
\psi(x) + c(x,y) \leq \psi(z) + c(z,y) \qquad \forall z \in \R^n.
$$
The functions $\psi$ and $\psi^c$ are said to be \textit{$c$-conjugate}.}
\end{definition}

Let us denote by $P_c(\R^n)$ the set of compactly supported probability measures. From now on, $\supp(\mu_0)$ and $\supp(\mu_1)$ will denote the supports of $\mu_0$ and $\mu_1$ respectively, i.e. the smallest closed sets on which $\mu_0$ and $\mu_1$ are respectively concentrated. Kantorovitch duality can be stated as follows (see \cite[Theorem 5.10]{villaniSF}):

\begin{theorem}
\label{THMduality} 
Let $\mu_0,\mu_1 \in P_c(\R^n)$ and $c:\R^n \times \R^n \rightarrow [0,\infty)$ be a lower continuous function. Then
there exists a $c$-convex function $\psi : \R^n \rightarrow \R$ such that the following holds: a transport plan
$\g \in \Pi(\mu_0,\mu_1)$ is optimal if and only if $\g(\p_c \psi)=1$
(that is, $\g$ is concentrated on the $c$-subdifferential of
$\psi$). Moreover $\psi$ can be chosen such that
\begin{eqnarray}\label{THMduality1}
\psi(x)&=&\sup_{y \in \supp(\mu_1)} \left\{ \psi^c(y) - c(x,y) \right\} \qquad
\forall x \in \R^n,
\\
\label{THMduality2}
\psi^c(y)&=&\inf_{x \in \supp(\mu_0)} \left\{\psi(x) + c(x,y)\right\} \qquad
\forall y \in \R^n.
\end{eqnarray}
\end{theorem}

By the above theorem we see that, in order to prove existence and
uniqueness of optimal transport maps, it suffices to prove that
there exist two Borel sets $Z_0,Z_1 \subset \R^n$, with $\mu_0(Z_0)=\mu_1(Z_1)=1$, such
that and $\p_c\psi$ is a graph inside $Z_0 \times Z_1$ (or
equivalently that $\p_c\psi(x) \cap Z_1$ is a singleton for all $x \in Z_0$).

\section{Proof of Theorem \ref{THM1}}\label{SPROOFTHM1}

First, we prove that (\ref{costTHM1}) holds. 
This is not new, but the formulas are not given in textbooks \textit{in extenso}.
The case $W=0$ is more classical and one may find, for instance in \cite{Lee-Mar67} the expression of $c$ in
that case:
\begin{equation}
\label{eq:costW0}
c(x,y)=\langle x-e^{-A}y\,,\,\mathcal{G}^{-1}(x-e^{-A}y)\rangle\ \ \mbox{with}\ \ \mathcal{G}=\int_0^1
e^{\tau A^*}B U^{-1} B^* e^{\tau A}\mathrm{d}\tau\ 
\end{equation}
where the controllability Grammian $\mathcal{G}$ is positive definite because the system is controllable.
This is indeed of the form (\ref{costTHM1}). Let us derive the general form from the classical linear
Pontryagin Maximum principle (see for instance the same reference).

The pseudo Hamiltonian $H_0 : \R^n \times \R^n \times \R^m \rightarrow \R$ associated with the optimization problem under study is given by 
\begin{eqnarray}
\label{pseudo}
H_0(x,p,u) := \langle p, Ax+Bu\rangle - \frac{1}{2} \langle x,Wx\rangle - \frac{1}{2} \langle u,Uu\rangle,
\end{eqnarray}
the control $u$ that maximizes this expression for each $(p,x)$ is given by
\begin{eqnarray}\label{contoptim}
u = U^{-1} B^*p
\end{eqnarray}
and the Hamiltonian $H:\R^n \times \R^n \rightarrow \R$, defined as $H(x,p)= \max \left\{ H_0(x,p,u) \, \vert \, u\in \R^m \right\}$ is given by 
\begin{eqnarray}
\label{ham}
H(x,p) 
& = &  \langle p, Ax\rangle - \frac{1}{2} \langle x,Wx\rangle + \frac{1}{2} \langle p, BU^{-1}B^*p\rangle.
\end{eqnarray}
Therefore the Hamiltonian differential equation $\dot x=\partial H/\partial p$, $\dot p=-\partial H/\partial
x$ associated to our minimization problem is given by 
\begin{eqnarray}
\label{hamsys}
\left(
\begin{array}{c}
  \dot{x} \\ \dot{p}
\end{array}
\right)
=
\left( \begin{array}{cc}
A & BU^{-1}B^* \\
W & -A^* 
\end{array} \right)
\left(
\begin{array}{c}
  x\\p
\end{array}
\right)
\end{eqnarray}
Denote by $R(\cdot): [0,1] \rightarrow M_{2n}(\R)$ the fundamental solution to the Cauchy problem 
\begin{eqnarray}
\label{cauchyfundamental}
\dot{R}(t) = \left( \begin{array}{cc}
A & BU^{-1}B^* \\
W & -A^* 
\end{array} \right) R(t) \quad \forall t \in [0,1], \qquad R(0)=I_{2n},
\end{eqnarray}
and write it as
$$
R(t) =  \left( \begin{array}{cc}
R_1(t) & R_2(t) \\
R_3(t) & R_4(t) 
\end{array} \right) \qquad \forall t \in [0,1],
$$
where $R_i(\cdot)$ is a $n \times n$ matrix for $i=1, \cdots 4$.  For every $x\in \R^n$ fixed, denote by $\exp_x$ the mapping from $\R^n$ to $\R^n$ which sends via the Hamiltonian vector field the initial adjoint vector $p\in \R^n$ to the final state $x(1)$, that is
\begin{equation}
\label{EXPx}
\exp_x (p) := R_1(1) x + R_2(1) p \qquad \forall p \in \R^n.
\end{equation}
Controllability of the system (\ref{system}) implies that the affine mapping $\exp_x$ is onto, and hence a bijection.
Indeed, there is, for any $x, y \in \R^n$, (at least) one control $u\in L^2\bigl([0,1];\R^m\bigr)$ which minimizes the cost 
$$
\int_{0}^{1} L( x(t;x,u),u(t)) dt 
$$
among all controls steering $x$ to $y$ in time $1$; thanks to the linear maximum principle, there corresponds
to each minimizing control a $p\in \R^n$ such that $\exp_x(p) =y$. Hence the matrix $R_2(t)$ is invertible for all $t$.

Given $x$ and $y$ in $\R^n$, the cost between $x$ and $y$ is given by
$$
c(x,y) = \int_0^1 \frac{1}{2} \langle x(t),Wx(t)\rangle +  \frac{1}{2} \langle p(t),BU^{-1}B^* p(t) \rangle dt,
$$
where $\bigl(x(\cdot),p(\cdot)\bigr) : [0,1] \rightarrow \R^n \times \R^n$ is the solution to the Hamiltonian system (\ref{hamsys}) satisfying
$$
x(0) = x \quad \mbox{ and } \quad \exp_x \bigl(p(0)\bigr) = y.
$$
Set $p:=p(0)$, that is
\begin{equation}
\label{VALp}
p= R_2(1)^{-1} \Bigl( y - R_1(1) x \Bigr).
\end{equation}
We deduce easily that $c$ is a smooth function of  the form
\begin{eqnarray}\label{formulacost}
c(x,y) = \frac{1}{2} \langle x,Q_1 x\rangle + \langle x,Cp\rangle + \frac{1}{2} \langle p, Q_2p\rangle \qquad \forall x,y \in \R^n, 
\end{eqnarray}
where both $Q_1,Q_2$ are non-negative symmetric $n\times n$ matrices given by
$$
Q_1 := \int_0^1 \Bigl( R_1(t)^* W R_1(t) + R_3(t)^* BU^{-1}B^* R_3(t)\Bigr) dt,
$$
$$
Q_2 :=  \int_0^1 \Bigl( R_2(t)^* W R_2(t) + R_4(t)^* BU^{-1}B^* R_4(t)\Bigr) dt,
$$
and where $C$ is the $n\times n$ matrix given by
$$
C :=  \int_0^1 \Bigl( R_1(t)^* W R_2(t) + R_3(t)^* BU^{-1}B^* R_4(t)\Bigr) dt.
$$
Let $x, y\in \R^n$ and $u \in L^2\bigl( [0,1]; \R^m\bigr)$ a control minimizing $c(x,y)$ be fixed.  For every control $v \in L^2\bigl( [0,1], \R^m\bigr)$, denote by $y(\cdot;y,v) : [0,1] \rightarrow \R^n$ the unique solution to the Cauchy problem
\begin{eqnarray*}
\left\{ 
\begin{array}{l}
\dot{y}(t) = A y(t) + B v(t) \quad \mbox{for a.e. } t \in [0,1],\\
y(1) =y.
\end{array}
\right.
\end{eqnarray*}
By definition of the cost $c$, there holds
\begin{eqnarray*}
\int_0^1  L( y(t;y,v),v(t)) dt - c\bigl( y(0;y,v),y\bigr) \geq 0 \qquad \forall v \in L^2\bigl( [0,1], \R^m\bigr).
\end{eqnarray*}
Moreover, there is equality in the above inequality whenever $v=u$. Thanks to the linear maximum principle, since the cost $c$ is smooth, this means that
$$
y = \exp_x(p) \quad \mbox{ with } \quad p = - \nabla_x c(x,y).
$$
Computing $\nabla_x c(x,y)$ by differentiating (\ref{formulacost})-(\ref{VALp}) with respect to the variable $x$ yields 
$$
Q_1x + C p - R_1(1)^* \bigl( R_2(1)^{-1}\bigr)^*\bigl(C^* x+ Q_2p\bigr)  = -p, 
$$
and finally (recall that $Q_1$ is symmetric)~:
\begin{eqnarray}\label{formulaC}
C=- I_n +  R_1(1)^* \bigl( R_2(1)^{-1}\bigr)^* Q_2, \quad Q_1= R_1(1)^* \bigl( R_2(1)^{-1}\bigr)^*C^*=CR_2(1)^{-1}R_1(1) .
\end{eqnarray}
Plugging this and (\ref{VALp}) into (\ref{formulacost}) yields the expression (\ref{costTHM1}) for $c(x,y)$
with
\begin{equation}
\label{DEF}
D= R_2(1)^{-1}R_1(1)\,,\ \ E=R_2(1)^{-1}\,,\ \ F=(R_2(1)^{-1})^*Q_2R_2(1)^{-1}\,,
\end{equation}
where $D$ is symmetric because (\ref{formulaC}) implies $Q_1=-E+E^*Q_2E$ and $Q_1,Q_2$ are symmetric; Positive definiteness of $D$
and $F$ may be deduced from the fact that there is no solution of (\ref{system}) driving 0 to a nonzero $y$ or a
nonzero $x$ to 0 in finite time with a zero cost.

Let us now show the existence and uniqueness of an optimal transport map. Let $\mu_0, \mu_1$ be two compactly supported probability measures in $\R^n$ and assume that $\mu_0$ is absolutely continuous with respect to the Lebesgue measure. Let $\psi, \phi:=\psi^c : \R^n \rightarrow \R$ be the Kantorovitch potentials given by Theorem \ref{THMduality} (note that $c$ is non-negative valued). First, since $c$ is smooth and both sets $\supp(\mu_0), \supp(\mu_1)$ are assumed to be compact, (\ref{THMduality1})-(\ref{THMduality2}) imply that  both  potentials $\psi, \phi$ are locally Lipschitz. Therefore, thanks to Rademacher's Theorem and the fact that $\mu_0$ is absolutely continuous with respect to the Lebesgue measure, $\psi$ is differentiable $\mu_0$-almost everywhere. Let $x \in  \supp(\mu_0)$ be a differentiability point of $\psi$. Let $y \in \supp(\mu_1)$ be such that
$$
\phi(y) - \psi (x) = c(x, y),
$$
and $u \in L^2\bigl( [0,1], \R^m\bigr)$ be a control minimizing 
$$
c(x,y) = \inf \left\{  \int_{0}^{1} L( x(t;x,v),v(t)) dt  \, \vert \, v\in L^2\bigl( [0,1];\R^m\bigr) \, \mbox{ s.t. } x(1;x,v) = y \right\}.
$$
Then we have 
\begin{eqnarray*} 
\int_0^1  L( y(t;y,v),v(t)) dt  \geq c\bigl(y(0;y,v),y\bigr) \geq      \phi(y)   - \psi \bigl(y(0;y,v)\bigr)   \qquad \forall v \in L^2\bigl( [0,1]; \R^m\bigr),
\end{eqnarray*}
with equality if $v=u$.  As above, this implies that 
$$
y = \exp_x(p) \quad \mbox{ with } \quad p = \nabla \psi (x).
$$
This shows that $y$ is uniquely determined for every differentiability point of $\phi$ and in turn proves existence and uniqueness of an optimal measurable transport map $x\mapsto y$. Note that the potential $\psi$ can be written as follows:
\begin{eqnarray*}
\psi (x) & = &  \sup_{y \in \supp(\mu_1)} \left\{ \phi(y) - c(x,y) \right\} \\
& = &  - \frac{1}{2} \langle x, D\, x\rangle +    \sup_{y \in \supp(\mu_1)} \left\{   \langle x, E\, y\rangle   - \frac{1}{2}  \langle y, F \,y\rangle  + \phi (y)\right\}.
\end{eqnarray*}
This shows that the function $\varphi: \R^n \rightarrow \R$ defined by
$$
\varphi ( x) :=  \psi(x) + \frac{1}{2} \langle x, D\, x\rangle   \qquad \forall x \in \R^n,
$$
is convex (as the \textit{sup} of a family of convex functions) while, for almost every $x\in
\supp(\mu_0)$, deriving the expression of $\exp_x(p)$ from (\ref{EXPx}) and (\ref{DEF}),
$$
T(x) = \exp_x \bigl(  \nabla \psi(x) \bigr) = E^{-1}\bigl( D x + \nabla \psi(x) \bigr) = E^{-1} \nabla\varphi(x)\;.
$$

\medskip

It remains to show that under additional assumptions, $T$ is continuous. 
One obviously has $\nabla\varphi(x)=E\;T(x)$ for all $x$; it is clear that this and $T_\sharp\mu_0=\mu_1$
imply $(\nabla\varphi)_\sharp\mu_0=\hat{\mu}_1$ for the measure $\hat{\mu}_1$ defined by
$$
\hat{\mu}_1 (A) = \mu_1 \Bigl( \left\{  E^{-1} x \, \vert \, x \in A \right\} \Bigr)
$$
for every Borel set $A \subset \R^n$ ($\hat{\mu}_1$ is the pushforward of $\mu_1$ by the map $x\mapsto Ex$).
Since $\varphi$ is convex, $\nabla\varphi$ is then the optimal transport map from $\mu_0$ to $\hat{\mu}_1$ for the cost
$\hat{c}(x,y) = \frac12|x-y|^2$. 
This implies the continuity of $\nabla\varphi$ ---hence of $T$--- according to
\cite[Theorem 12.50 (i)]{villaniSF} which asserts the following: 
if $\hat{\mu}_0, \hat{\mu}_1$ are two compactly supported probability measures in $\R^n$ 
associated with densities $\hat{f}_0, \hat{f}_1$ which are bounded from below and above 
on $\supp(\hat{f}_0)$ and $\supp(\hat{f}_1)$ respectively, 
and if $\supp(\hat{f}_0)$ is connected and $\supp(\hat{f}_1)$ is convex, 
then the optimal transport map $\hat{T}: \R^n \rightarrow \R^n$ from $\hat{\mu}_0$ to 
with respect to the  Euclidean quadratic cost $\hat{c}(x,y)= \frac 1 2 |x-y|^2$ is continuous. 
Take $\hat{\mu}_1$ defined above and set $\hat{\mu}_0=\mu_0$; they do satisfy these assumptions for $\mu_0$
and $\mu_1$ and applying an invertible linear map does not change convexity of
the support of a density, and we just proved that $\hat{T}$ is $\nabla\varphi$ with the same $\varphi$ as above.



\section{Proof of Theorem \ref{THM2}}\label{SPROOFTHM2}

Let us first treat the case $d=0$ separately. According to (\ref{eq:kalman}), either there is no control
($m=0$) or the matrix $B$ is zero; in both cases the system reads $\dot{x}=Ax$ and
\begin{eqnarray*}
c (x,y) = \left\{ \begin{array}{cl}
\frac{1}{2} \int_0^1 \left\langle e^{tA} \, x, W e^{tA} \, x \right\rangle dt & \mbox{ if } y = e^{A}x\,,  \\
+ \infty & \mbox{ otherwise}.
\end{array}
\right.
\end{eqnarray*}
Then there is only one map $S$ such that $\mbox{cost}_c (S) <\infty$, given by
$S(x) = e^A x$ for all $x$, the compatibility on measures is $\mu_1=S_\sharp\mu_0$ with this precise $S$, and the result is proved in the case $d=0$. 

From now on, we assume that $d\in \{1,\cdots, n-1\}$,
i.e. none of the two blocks in (\ref{system2}) is void ($d\geq1$, $n-d\geq1$).
For every $x=\bigl( x_1,x_2\bigr) \in \R^d \times \R^{n-d} = \R^n$ and $u\in L^2 \bigl( [0,1];\R^m\bigr),$ the
solution $x(\cdot;x,u) :[0,1] \rightarrow \R^n$ to the Cauchy problem (\ref{cauchy}) with the decomposition
given by (\ref{system2}) satisfies 
\begin{eqnarray}\label{solutionTHM2}
x(t;x,u) = \left( \begin{matrix}
x_1 (t;x,u ) \\
x_2 \bigl(t;x_2\bigr)
\end{matrix}
\right) = \left( \begin{matrix}
e^{tA_1} x_1 + e^{tA_1} \int_0^t e^{-sA_1} \left[ A_3 e^{sA_2} x_2 + B_1 u(s)\right] ds  \\
e^{tA_2} x_2 
\end{matrix}
\right), 
\end{eqnarray}
for every $t\in [0,1]$.  Denote by $\bar{x}_1 (\cdot)  = \bar{x}_1 \bigl( \cdot ; x_1, u\bigr) : [0,1] \rightarrow \R^d$ the solution to the Cauchy problem
\begin{eqnarray}\label{cauchyTHM2}
\left\{ 
\begin{array}{l}
\dot{\bar{x}}_1(t) = A_1 \bar{x}_1(t) + B_1 u(t) \quad \mbox{for a.e. } t \in [0,1],\\
\bar{x}_1(0) = x_1.
\end{array}
\right.
\end{eqnarray}
Then we have 
$$
x_1 \bigl(t;x,u\bigr) = \bar{x}_1 \bigl( t;x_1,u\bigr) +  G(t) x_2 \qquad \forall t \in [0,1],
$$
where $G:[0,1] \rightarrow M_{d,n-d}(\R)$ is defined by
$$
G (t) := e^{tA_1} \int_0^t e^{-sA_1} A_3 e^{sA_2} ds \qquad \forall t \in [0,1].
$$

Therefore we have for almost every $t\in [0,1]$,
\begin{eqnarray*}
L \bigl( x(t;x,u),u(t) \bigr)
& = & \frac{1}{2} \left\langle  x(t;x,u), W  x(t;x,u) \right\rangle + \frac{1}{2} \langle u(t), U u(t) \rangle \\
& = &  \frac{1}{2} \left\langle \bar{x}_1 \bigl(t;x_1,u\bigr), W_1  \, \bar{x}_1\bigl(t;x_1,u\bigr) \right\rangle +  \frac{1}{2} \langle u(t), U u(t) \rangle \\
& \quad & \qquad  \qquad \qquad \qquad  + \left\langle  \bar{x}_1 \bigl(t;x_1,u\bigr), X \bigl(t;x_2\bigr) \right\rangle + l \bigl(t;x_2\bigr),
\end{eqnarray*}
where the symmetric matrix $W$ is decomposed, in the coordinates of (\ref{system2}) as
$$W=\left(
\begin{array}{cc}
  W_1&W_3\\{W_3}^*&W_2
\end{array}
\right)$$
with $W_1$ and $W_2$ positive definite $d\times d$ and $n-d\times n-d$ matrices respectively, $W_3$ a $d\times
n-d$ matrix, and
\begin{eqnarray*}
X\bigl(t;x_2\bigr)&=&W_1\, G(t) x_2 +W_3 \,x_2 \bigl(t;x_2 \bigr)\,,
\\
l\bigl(t;x_2\bigr) &=& \frac{1}{2} \left\langle 
  \left(  \begin{matrix} G(t) x_2 \\ x_2 \bigl(t;x_2 \bigr)\end{matrix} \right)
   , 
  \left( \begin{array}{cc} W_1&W_3\\{W_3}^*&W_2 \end{array} \right)
  \left(  \begin{matrix} G(t) x_2 \\ x_2 \bigl(t;x_2\bigr)  \end{matrix} \right)  
  \right\rangle\,.
\end{eqnarray*}
In conclusion, we have for every $x=\bigl( x_1, x_2\bigr), y=\bigl( y_1, y_2\bigr) \in \R^n$,
\begin{eqnarray}\label{costTHM2}
c (x,y) = \left\{ \begin{array}{cl}
\bar{c}_{x_2} \bigl( x_1, y_1- G(1)x_2\bigr) + \int_0^1 l \bigl(t;x_2\bigr)dt & \mbox{ if } y_2 = e^{A_2} x_2 \\
+ \infty & \mbox{ otherwise},
\end{array}
\right.
\end{eqnarray}
where the cost $\bar{c}_{x_2} : \R^d \times \R^d \rightarrow [0, +\infty)$ is defined by 
\begin{multline}\label{costLQTHM2}
\bar{c}_{x_2}   \bigl( x_1, z_1\bigr) := \\
\inf \left\{  \int_{0}^{1} \bar{L}_{x_2} \bigl(t, \bar{x}_1(t;x_1,u),u(t) \bigr) dt  \, \vert \, u\in L^2\bigl( [0,1];\R^m\bigr) \, \mbox{ s.t. } \bar{x}_1 \bigl(1;x_1,u\bigr) = z_1 \right\},
\end{multline}
for any $x_1, z_1 \in \R^d,$ and where the Lagrangian $\bar{L}_{x_2} : \R \times \R^d \times \R^m \rightarrow [0,+\infty)$ is defined by
\begin{eqnarray}
\bar{L}_{x_2} \bigl(t, z,u\bigr) :=  \frac{1}{2} \left\langle z, W_1  \, z \right\rangle +  \frac{1}{2} \left\langle u(t), U u(t) \right\rangle + \left\langle z, X \bigl(t;x_2\bigr)  \right\rangle.
\end{eqnarray}
We proceed now as in the proof of Theorem \ref{THM1}. Let $x_2 \in \R^{n-d}$ be fixed, the pseudo Hamiltonian $H_0 : \R \times \R^d \times \R^d \times \R^m \rightarrow \R$ associated with the cost $\bar{c}_{x_2}$ is given by 
\begin{eqnarray}
\label{pseudo2}
H_0(t,z,p,u) := \left\langle p, A_1 z+ B_1 u \right\rangle - \frac{1}{2} \left\langle z,W_1 \, z \right\rangle - \frac{1}{2} \left\langle u,U \, u \right\rangle - \left\langle z, X\bigl(t;x_2\bigr) \right\rangle.
\end{eqnarray}
Then $\frac{\partial H_0}{\partial u}=0$ yields $u = U^{-1} B_1^*p$ and the Hamiltonian $H: \R \times \R^d \times \R^d \rightarrow \R$ is given by 
\begin{eqnarray*}
H(t,z,p) & = & \max \left\{ H_0(t,z,p,u) \, \vert \, u\in \R^m \right\} \\
& = &  \left\langle p, A_1 \, z \right\rangle - \frac{1}{2} \left\langle z, W_1 \, z \right\rangle + \frac{1}{2} \left\langle p, B_1 U^{-1}B_1^*p \right\rangle -  \left\langle z, X\bigl(t;x_2\bigr) \right\rangle.
\end{eqnarray*}
Therefore the Hamiltonian system associated to our minimization problem is given by 
\begin{eqnarray}
\label{hamsysTHM2}
\left\{ 
\begin{array}{lll}
\dot{z} & = & A_1 z + B_1U^{-1}B_1^*p \\
\dot{p} & = & - A_1^*p +W_1 z + X \bigl( t;x_2\bigr).
\end{array}
\right.
\end{eqnarray}
Denote by $R(\cdot): [0,1] \rightarrow M_{2d}(\R)$ the fundamental solution to the Cauchy problem 
\begin{eqnarray}
\label{cauchyfundamental2}
\dot{R}(t) = \left( \begin{array}{cc}
A_1 & B_1U^{-1}B_1^* \\
W_1 & -A_1^* 
\end{array} \right) R(t) \quad \forall t \in [0,1], \qquad R(0)=I_{2d},
\end{eqnarray}
and write it as
$$
R(t) =  \left( \begin{array}{cc}
R_1(t) & R_2(t) \\
R_3(t) & R_4(t) 
\end{array} \right) \qquad \forall t \in [0,1],
$$
where $R_i(\cdot)$ is a $d \times d$ matrix for $i=1, \cdots 4$. Any solution $\bigl( z(\cdot),p(\cdot)\bigr) : [0,1] \rightarrow \R^d \times \R^d$ of (\ref{hamsysTHM2}) with $z(0)= z, p(0)=p$ can be written as
$$
\left\{
\begin{array}{lll}
z(t) & = & R_1(t) z + R_2(t) p + \bar{z}_{x_2}(t) \\
p(t) &= & R_3(t) z + R_4(t) p + \bar{p}_{x_2}(t)
\end{array}
\right.
\qquad \forall t\in [0,1],
$$
where $\bigl( \bar{z}_{x_2}(\cdot),\bar{p}_{x_2}(\cdot)\bigr) : [0,1] \rightarrow \R^d \times \R^d$ of (\ref{hamsysTHM2}) satisfying $\bar{z}_{x_2}(0)=0,  \bar{p}_{x_2}(0)=0$. For every $x_1 \in \R^d$, the mapping $\exp_{x_1} : \R^d \rightarrow \R^d$ defined by 
$$
\exp_{x_1} (p) := R_1(1) x_1 + R_2(1) p + \bar{z}_{x_2} (1) \qquad \forall p \in \R^d,
$$
is an affine bijection. For every  $x_1, z_1, p \in \R^d$ with $z_1 = \exp_{x_1} (p)$, there holds 
\begin{eqnarray}\label{formulacostTHM2}
\bar{c}_{x_2} \bigl(x_1,z_1\bigr) = \frac{1}{2} \left\langle x_1 ,Q_1 \, x_1 \right\rangle +  \frac{1}{2} \left\langle p, Q_2 \, p \right\rangle + \left\langle x_1, C \, p \right\rangle +  \left\langle x_1, v_{x_2} \right\rangle + \left\langle p, w_{x_2} \right\rangle + k_{x_2},
\end{eqnarray}
where both $Q_1,Q_2$ are non-negative symmetric $d\times d$ matrices given by
$$
Q_1 := \int_0^1 \Bigl( R_1(t)^* W_1 R_1(t) + R_3(t)^* B_1U^{-1}B_1^* R_3(t)\Bigr) dt,
$$
$$
Q_2 :=  \int_0^1 \Bigl( R_2(t)^* W_1 R_2(t) + R_4(t)^* B_1U^{-1}B_1^* R_4(t)\Bigr) dt,
$$
where $C$ is the $d \times d$ matrix given by
$$
C :=  \int_0^1 \Bigl( R_1(t)^* W_1 R_2(t) + R_3(t)^* B_1 U^{-1}B_1^* R_4(t)\Bigr) dt,
$$
and the vectors $v_{x_2}, w_{x_2} $ are given by
$$
v_{x_2} := \int_0^1 \Bigl( R_1(t)^* W_1 \, \bar{z}_{x_2}(t) + R_3(t)^* B_1U^{-1}B_1^* \, \bar{p}_{x_2}(t)  +  R_1(t)^* \, X \bigl( t;x_2 \bigr) \Bigr)dt,
$$
$$
w_{x_2} := \int_0^1 \Bigl( R_2(t)^* W_1 \, \bar{z}_{x_2}(t) + R_4(t)^* B_1U^{-1}B_1^* \, \bar{p}_{x_2}(t)+ R_2(t)^* \, X \bigl( t;x_2 \bigr)   \Bigr) dt,
$$
and 
$$
k_{x_2} := \int_0^1 \left( \frac{1}{2} \left\langle \bar{z}_{x_2}(t), W_1 \, \bar{z}_{x_2}(t) \right\rangle  + \frac{1}{2} \left\langle \bar{p}_{x_2}(t), B_1U^{-1}B_1^* \,  \bar{p}_{x_2} (t) \right\rangle +  \left\langle \bar{z}_{x_2}(t), X \bigl( t;x_2 \bigr)  \right\rangle \right) dt.
$$
We proceed now as in the proof of Theorem \ref{THM1}. Since $\bar{c}_{x_2}$ is smooth, using the linear maximum principle, taking the derivative in the $x_1$ variable in (\ref{formulacostTHM2}) and using that 
$p = R_2(1)^{-1} \bigl[ z_1 - R_1(1)x_1 - \bar{z}_{x_2}(1) \bigr]$ yields
\begin{displaymath}
Q_1 \, x_1 + C \,p +v_{x_2}
- R_1(1)^* \bigl( R_2(1)^{-1}\bigr)^* \left(Q_2 \, p + C^*\,x_1+ w_{x_2}\right)=-p
\end{displaymath}
and finally $Q_1-R_1(1)^*(R_2(1)^{-1})^*C^*=C-R_1(1)^*(R_2(1)^{-1})^*Q_2+I=v_{x_2}- R_1(1)^* w_{x_2}=0$, whence
(recall that $Q_2$ is symmetric):
\begin{eqnarray}\label{formulaC2}
v_{x_2} = R_1(1)^* \bigl( R_2(1)^{-1}\bigr)^* \, w_{x_2}, \quad  \left\{ \begin{array}{l}
C=- I_n +  R_1(1)^* \bigl( R_2(1)^{-1}\bigr)^* Q_2 \\
Q_1=R_1(1)^*(R_2(1)^{-1})^*C^*=CR_2(1)^{-1}R_1(1)
\end{array}
\right.
\end{eqnarray}
In conclusion, setting $D := R_2 (1)^{-1}R_1(1), E :=R_2 (1)^{-1}$, and $F:= E^*Q_2E$ yields
\begin{multline*}
\bar{c}_{x_2} \bigl(x_1,z_1\bigr) =  \frac{1}{2} \left\langle x_1, D \, x_1 \right\rangle -  \left\langle x_1, E \,  \bigl( z_1 -\bar{z}_{x_2}(1)\bigr) \right\rangle \\
+  \frac{1}{2} \left\langle  \bigl( z_1 -\bar{z}_{x_2}(1)\bigr), F \,  \bigl( z_1 -\bar{z}_{x_2}(1)\bigr) \right\rangle + \left\langle E  \bigl( z_1 -\bar{z}_{x_2}(1)\bigr), w_{x_2} \right\rangle + k_{x_2},
\end{multline*}
which in turn implies (by (\ref{costTHM2}))
\begin{multline}\label{costfinalTHM2}
c \left( \bigl( x_1,x_2\bigr), \bigl( y_1, e^{A_2}x_2\bigr) \right) = \\
\frac{1}{2} \left\langle x_1, D \, x_1 \right\rangle -  \left\langle x_1, E \,  \bigl( y_1 - G(1) x_2 -\bar{z}_{x_2}(1)\bigr) \right\rangle \\
+  \frac{1}{2} \left\langle  \bigl( y_1 -G(1) x_2 -\bar{z}_{x_2}(1)\bigr), F \,  \bigl( y_1 - G(1) x_2 -\bar{z}_{x_2}(1)\bigr) \right\rangle \\
+ \left\langle E  \bigl( y_1 -G(1) x_2 -\bar{z}_{x_2}(1)\bigr), w_{x_2} \right\rangle + k_{x_2} + \int_0^1 l\bigl(t;x_2\bigr) dt.
\end{multline}
$D$ and $F$ are symmetric definite positive for the same reason as in the proof of Theorem~\ref{THM1}.

We are now ready to prove the result. 
For every $\bar{x}_2 \in \R^{n-d}$, set
\begin{eqnarray*}
M_{\bar{x}_2} &:=& \left\{ x=(x_1,x_2) \in \R^d \times \R^{n-d} \, \vert \, x_2 = \bar{x}_2\right\},
\\
N_{\bar{x}_2} &:=& \left\{ x=(x_1,x_2) \in \R^d \times \R^{n-d} \, \vert \, x_2 = e^{A_2} \bar{x}_2\right\},
\end{eqnarray*}
and let $\mu_0= f_0 \mathcal{L}^n, \mu_1= f_1 \mathcal{L}^n$ with densities $f_0, f_1 \in L^1(\R^n)$ be two
compactly supported probability measures. 

\smallskip

\noindent\textbf{Fact:} {\it A measurable map $S$ satisfies $S_\sharp\mu_0=\mu_1$ and $\mbox{cost}_c (S) < +
\infty$  if and only if there exists, for almost all $x_2$, a measurable $S_{x_2}:\R^d\to\R^d$ such that 
\begin{equation}
\label{eq:jb1}
S (x_1,x_2) = \Bigl( S_{x_2} (x_1), e^{A_2} x_2 \Bigr)  \qquad \mbox{ for a.e. } x_1
\end{equation}
and $S_{x_2}$ pushes forward the measure $\mu_0^{x_2}$ on $M_{x_2}$ defined by
$$
\mu_0^{x_2} := f_0 \bigl(\cdot, x_2\bigr) \left| \det \bigl(e^{-A_2}\bigr) \right| \, dx_1
$$
to the measure $\mu_1^{x_2}$ on $N_{x_2}$ defined by
$$
\mu_1^{x_2} := f_1 \left(\cdot, e^{A_2} x_2\right) \, dx_1.
$$
}
Indeed,
by Fubini's Theorem, the map $x_1\mapsto S(x_1,x_2)$ is measurable for almost all $x_2$ and 
$\mbox{cost}_c (S) < +\infty$ implies
$$
\int_{\R^d} c\bigl( (x_1,x_2), S(x_1,x_2) \bigr) f_0(x_1,x_2) dx_1 < + \infty \qquad \mbox{ for a.e. } x_2
$$
and in turn, from (\ref{costTHM2}), this implies, for almost every $x_2 \in \R^{n-d}$, that
$S(x_1,x_2) \in N_{x_2}$ for almost all $x_1$,
which means that $S(x_1,x_2)$ has the form (\ref{eq:jb1}).
Now $S_\sharp \mu_0=\mu_1$ means
$
\int_{R^n} h(x) d\mu_1(x) = \int_{R^n} h\bigl(S(x)\bigr) d\mu_0(x)
$
for any $h \in L^1\bigl(\R^n\bigr)$.
From Fubini's Theorem together with a change of variable $y_2=e^{A_2}x_2$ in the left-hand side, this yields, for any positive $h \in L^1\bigl(\R^n\bigr)$
\begin{eqnarray*}
& \quad &   \int_{\R^{n-d}} \left( \int_{\R^d} h(x_1,x_2) f_1(x_1,x_2) dx_1 \right) dx_2 \\
& = &  \int_{\R^{n-d}} \left( \int_{\R^d} h\left(   \Bigl( S_{e^{-A_2} y_2} (x_1), y_2 \Bigr) \right)  f_0\left(x_1, e^{-A_2}y_2 \right) \left| \det \bigl(e^{-A_2}\bigr) \right| dx_1 \right) dy_2\;.
\end{eqnarray*}
This is equivalent to $(S_{x_2})_\sharp \mu_0^{x_2}=\mu_1^{x_2}$ and proves the fact.

\smallskip
If $S$ satisfies $S_\sharp\mu_0=\mu_1$ and $\mbox{cost}_c (S) < +\infty$, the form (\ref{eq:jb1}) implies
that, for all measurable $h$ that depends on $x_2$ only, one has 
$$\int_{\R^d\times\R^{n-d}}h(x_2)d\mu_1(x_1,x_2)=\int_{\R^d\times\R^{n-d}}h(e^{A_2}x_2)d\mu_0(x_1,x_2),$$
and this implies (\ref{HYPTHM2bis}). Conversely, assume now that the measures satisfy (\ref{HYPTHM2bis}); since
$$
x_2\mapsto \int_{M_{x_2}}d\mu^{x_2}_0\quad\mbox{and}\quad
x_2\mapsto \int_{N_{x_2}}d\mu^{x_2}_1
$$
are the densities of $(e^{A_2}\circ\pi_2)_\sharp \mu_0$ and $(\pi_2)_\sharp \mu_1$ respectively,
the measures $\mu^{x_2}_0$ and $\mu^{x_2}_1$ have, for almost all $x_2$, the same total mass,
hence can be considered as probability measures after normalisation.
Then using (\ref{costfinalTHM2}) and arguing as in the proof of Theorem \ref{THM1}, we deduce that for almost
every $x_2 \in \R^{n-d}$, there is a unique optimal transport map $T_{x_2} : M_{x_2} \rightarrow N_{x_2}$ together with a convex function $\varphi_{x_2} : \R^d \rightarrow \R$ such that
$$
T_{x_2} \bigl(x_1\bigr) =  \Bigl( E^{-1} \nabla \varphi_{x_2} (x_1), e^{A_2} x_2 \Bigr) \qquad \mbox{ for a.e. } x_1 \in \R^d.
$$
This shows the existence and uniqueness of an optimal transport map $T$, defined by $T(x_1,x_2)=(E^{-1} \nabla
\varphi_{x_2} (x_1), e^{A_2} x_2)$ from $\mu_0$ to $\mu_1$ (and a fortiori shows (\ref{HYPTHM2})).
Note that joint measurability of $T$ with respect to the two variables may be seen as a consequence of the necessity part of \cite[Theorem 3.2]{ap03} which asserts that any optimal plan is concentrated on a $c$-monotone Borel subset of $\R^n \times \R^n$.

If we assume in addition that $f_0, f_1$ are both continuous and bounded from below and above on their support, then each $T_{x_2}$ is continuous. This is a consequence of Theorem $1$ applied for each $x_2$ to  the mass transportation problem from $\mu_0^{x_2}$ to $\mu_1^{x_2}$ with respect to the cost given by (\ref{costfinalTHM2}). But by \cite[Corollary 5.23]{villaniSF}, there is stability of transport maps. If $\{x_2^k\}$ is a sequence converging to $x_2$, the measures $\mu_0^{x_2^k}$ and $\mu_1^{x_2^k}$ weakly converge to $\mu_0^{x_2}$ and $\mu_1^{x_2}$ respectively, and if in addition the corresponding costs converge uniformly as well, then the transport $\{T_{x_2^k}\}$ maps converge to $T_{x_2}$ in probability. Since all $T_{x_2}^k$ are indeed uniformly H\"older continuous (see \cite[Theorem 50 (i)]{villani03} or \cite[Theorem 12.50]{villaniSF}), and references therein), this concludes the proof of Theorem \ref{THM2}.

\section{Examples}\label{SEXAMPLES}

\subsection{The case $W=0$}

If $W=0$ then the second line in (\ref{hamsys}) is an independent first order differential equation in the $p$ variable; hence $R_3(t)=0$ for any $t\in [0,1]$ and the form (\ref{eq:costW0}) of the cost $c$.
If in addition, we assume that $A=0$ and that the system is controllable, then the matrix $B$ is necessarily invertible. In that case, we leave the reader to show that the $2n \times 2n$ matrix $R(t)$ has the form
$$
R(t) =  \left( \begin{array}{cc}
I_n & BU^{-1}B^* \\
0_n & I_n 
\end{array} \right). 
$$
Then there holds 
$$
c(x,y)  = \frac{1}{2}\left\langle \bigl(B U^{-1}B^*\bigr)^{-1}  \Bigl( y - x \Bigr),  \Bigl( y - x \Bigr)\right\rangle \qquad \forall x,y \in \R^n.
$$
And any optimal transport map given by Theorem \ref{THM1} has the form
$$
T = \bigl( B U^{-1}B^* \bigr) \, \nabla \varphi,
$$
where we used (\ref{DEF}), with $\varphi$ is a convex function.
This can also be viewed as a consequence of \cite{brenier91} because the above cost is the Euclidean norm corresponding to the positive matrix $(B U^{-1}B^*)^{-1}$ and is $\bigl( B U^{-1}B^* \bigr) \, \nabla \varphi$ is nothing but the gradient of $\varphi$ for this Euclidean metric.

\subsection{The case $A=0,\quad B=U=I_n$}
To have a nice form of $R(t)$ let us assume that $W$ is symmetric definite positive. In this case, we leave  the reader to show that the $2n \times 2n$ matrix $R(t)$ has the form :
$$
R(t) =  \left( \begin{array}{cc}
cosh\,\,\bigl(tW^{\frac{1}{2}}\bigl) & sinh\,\,\bigl(tW^{\frac{1}{2}}\bigl)W^{-\frac{1}{2}} \\
W^{\frac{1}{2}}.sinh\,\,\bigl(tW^{\frac{1}{2}}\bigl) &  cosh\,\,\bigl(tW^{\frac{1}{2}}\bigl)
\end{array} \right). 
$$
Where 
$$
cosh\,\,\bigl(tW^{\frac{1}{2}}\bigl)=\sum_{n=0}^\infty \frac{t^{2n}}{2n!}W^{n}, \quad\quad sinh\,\,\bigl(tW^{\frac{1}{2}}\bigl)= \sum_{n=0}^\infty \frac{t^{2n+1}}{2n+1\,!}W^{n+\frac{1}{2}}
$$
And any optimal transport map, given by Theorem \ref{THM1},  has the form
$$
T = \bigl( sinh\,\,W^{\frac{1}{2}}.W^{\frac{-1}{2}} \bigr) \, \nabla \varphi,
$$
where we used (\ref{DEF}) and that $\varphi$ is a convex function.

\section{Final comments}\label{Scomments}

It is worth noticing that the existence and uniqueness part in Theorem \ref{THM1} is not new. It is a consequence of \cite[Theorem 4.1]{al07} together with the Lipschitz regularity of the cost. The formula (\ref{costTHM1}) implies that the Ma-Trudinger-Wang tensor associated with the cost is identically zero. Such a result has been obtained previously by McCann and Lee  without computing explicitly the cost (see \cite[Theorem 1.1]{lm10}). We finally observe that (\ref{TransportTHM1}) means that the optimal transport map $T$ is convex up to a linear change of coordinates. Such a property is related to the vanishing of the Ma-Trudinger-Wang tensor that we mentioned above and \cite[Theorem 4.3]{fkm09}. We refer the interested reader to \cite{frv,villaniSF} for a more details about the Ma-Trudinger-Wang tensor and its link with regularity of optimal transport maps.


\end{document}